\documentclass[12pt]{article}
\usepackage{sectsty}
\usepackage{mathptmx}
\usepackage{fancyhdr, amssymb, amsmath, amsthm, enumerate, enumitem, lastpage}     
\usepackage{tikz, tikz-cd}
\usepackage{array}
\usepackage{multirow}
\usepackage{mathtools}
\usepackage{hyperref}
\usepackage{physics}
\usepackage{cite, url}
\usepackage{etoolbox} 
\preto\align{\par\nobreak\small\noindent} 
\expandafter\preto\csname align*\endcsname{\par\nobreak\small\noindent} 
\newcommand{\floor}[1]{\lfloor #1 \rfloor}

\newcommand{\ehk}[1]{e_{\mathrm{HK}}\left(#1\right)}
\newcommand{\ehkideal}[2]{e_{\mathrm{HK}}\left(#1, #2\right)}
\newcommand{\len}[2]{\lambda_{#1}\left(#2\right)}
\newcommand{\height}[1]{\mathrm{ht} \left(#1\right)}

\newcommand{\spec}[1]{\mathrm{Spec}\left(#1\right)}

\newcommand{\mf}[1]{\mathfrak{#1}}

\newcommand{\fbp}[1]{\left[ #1\right]}
\newcommand{\unideal}[1]{\left(\underline{#1}\right)}

\newcommand{\ol}[1]{\overline{#1}}
\newcommand{\indt}{\hspace{.3in}}

\newcommand{\ul}[1]{\underline{#1}}

\usepackage{graphicx}%
\usepackage{multirow}%
\usepackage{amsmath,amssymb,amsfonts}%
\usepackage{amsthm}%
\usepackage{mathrsfs}%
\usepackage[title]{appendix}%
\usepackage{xcolor}%
\usepackage{textcomp}%
\usepackage{manyfoot}%
\usepackage{booktabs}%
\usepackage{algorithm}%
\usepackage{algorithmicx}%
\usepackage{algpseudocode}%
\usepackage{listings}%
\usepackage{fullpage,xcolor,hyperref,verbatim}
\usepackage[sfdefault]{biolinum}
\usepackage{mathrsfs, multicol}
\usepackage{csquotes}
\newcommand{\exideal}[3]{\left(#1, \, #2\right)#3}
\newcommand{\Psifun}[3]{\varphi_{#1}\left(#2; \, #3\right)}
\newcommand{\PsifunZ}[1]{\Psifun{J}{R}{z^{#1}}}

\newtheorem{thm}{Theorem}[section]
\newtheorem*{thm*}{Theorem}
\newtheorem*{ded*}{Dedication}
\newenvironment{remark}[1][]{ \refstepcounter{thm}\par\medskip
   \noindent \textbf{Remark \arabic{section}.\arabic{thm}: #1} \rmfamily}{}

\newtheorem{lemma}[thm]{Lemma}
\newtheorem{conjecture}[thm]{Conjecture}
\newtheorem{prop}[thm]{Proposition}
\newtheorem{defn}[thm]{Defintion}
\newtheorem{cor}[thm]{Corollary}

\fancypagestyle{plain}{
  \fancyhf{}
  
  \fancyhead[R]{}
  \fancyfoot[R]{\thepage\  / \pageref{LastPage}}
}
\usepackage{etoolbox} 
\preto\align{\par\nobreak\small\noindent} 
\expandafter\preto\csname align*\endcsname{\par\nobreak\small\noindent} 
\begin{document}

\begin{center}
\text{\Large Bounds for the Hilbert-Kunz Multiplicity of Singular Rings} \\
\vspace{.5in}
\small{\textbf{Nicholas O. Cox-Steib} \,\, \textit{noc552@utulsa.edu} \\ University of Tulsa \\
\textbf{Ian M. Aberbach} \, \, \textit{aberbachi@missouri.edu} \\ University of Missouri, Columbia}
\end{center}

\vspace{.5in}

\begin{abstract}
	In this paper we prove that the Watanabe-Yoshida conjecture holds up to dimension $7$.  Our primary new tool is a function, $\Psifun{J}{R}{z^{t}},$ that interpolates between the Hilbert-Kunz multiplicities of a base ring, $R$, and various radical extensions, $R_n$.  We prove that this function is concave and show that it's rate of growth is related to the size of $\ehk{R}$.  We combine techniques from \cite{CHDZ} and \cite{IanNewEst} to get effective lower bounds for $\varphi,$ which translate to improved bounds on the size of Hilbert-Kunz multiplicities of singular rings.  The improved inequalities are powerful enough to show that the conjecture of Watanabe and Yoshida holds in dimension $7$.      \\
\end{abstract}

\vspace{.25in}
\begin{center}
	{\it This paper is dedicated to Ngo Viet Trung on the occasion of his 70th birthday.}
\end{center}
\vspace{.25in}

\\
\newpage 
\section{Introduction}
\subsection{The Watanabe-Yoshida Conjecture}
\indt	Let $(R, \mf{m}_R, \kappa)$ denote a local Noetherian ring of positive characteristic $p>0$ and dimension $d.$  The Hilbert-Kunz multiplicity of an $\mf{m}_R$-primary ideal, $J$, is given by the limit\footnote{This limit was shown to be well defined by Paul Monsky in \cite{monsky1983hilbert}.},
\begin{align*}
	\ehkideal{J}{R} := \lim_{e \to \infty} \dfrac{1}{p^{ed}} \len{R}{\dfrac{R}{J^{\fbp{p^e}}R}},
\end{align*}
where, $J^{\fbp{p^e}} = (x^{p^e} \, | \, x \in J),$ denotes the Frobenius bracket power of $J$ and we write $\lambda_R$ to denote length as an $R$-module. \\

The number, $\ehkideal{\mf{m}_R}{R} = \ehk{R},$ is an asymptotic measure of the flatness of the Frobenius map, $R \to R$, and thus gives a numerical measure of the 'singular-ness' of $R$.  It is straightforward to show that $\ehk{R} \ge 1,$ and, a Theorem of \cite{watanabe2000hilbert} shows that, when $R$ formally unmixed, it is regular if and only if $\ehk{R} = 1.$\footnote{The formally unmixed condition is necessary here; it is not difficult to construct examples of singular rings with embedded prime ideals that satisfy $\ehk{R} = 1.$  Recall that $(R, \mf{m}_R)$ is {\it formally unmixed} if the $\mf{m}_R$-adic completion, $\widehat{R},$ is unmixed --- that is, $\dim \widehat{R}/\mf{p} = \dim \widehat{R}$ for all associated primes of $\widehat{R}.$}   Heuristically one then expects that `the closer $\ehk{R}$ is to $1$, the closer the Frobenius morphism $R \to  R$ is to being flat,' or, equivalently, 'the larger the value of $\ehk{R}-1$ the more severe the singularity of $R$.'\\

It is therefore natural to ask what can be said about singularities, $(R, \mf{m}_R)$, that have $\ehk{R}$ close to one.  Such considerations led Watanabe and Yoshida to formulate the following conjecture, which appears as Conjecture 4.2 of \cite{WY3Dim}: \\

\begin{conjecture} \label{WY Conjecture} ({\bfseries The Watanabe-Yoshida Conjecture}) \\

	 Suppose that $(R, \mf{m}_R, \kappa)$ is a Noetherian local ring of dimension $d>0$ and characteristic $p>2.$  Assume $\kappa = \ol{\kappa},$ and that $R$ is formally unmixed.  Let $R_{p, d}$ denote the quadratic hypersurface
	\begin{align*}
		R_{p,d} = \kappa[[x_0, \dots, x_d]]/(x_0^2+\dots+x_d^2).	
	\end{align*}
If $R$ is not regular, then
	\begin{align} \label{WYineq}
		\ehk{R} \ge \ehk{R_{p,d}} \ge 1 + m_d,
	\end{align}
where the first inequality, $\ehk{R} \ge \ehk{R_{p,d}},$ is strict unless the $\mf{m}_R$-adic completion of $R$ is isomorphic to $R_{p,d}.$

\end{conjecture}

\vspace{.25in}

Here, the numbers $m_d$, which are defined by the relation 
\begin{align*}
	\sec x + \tan x = 1 + \sum_{d=1}^{\infty} m_d x^d \hspace{.5in} (-\pi/2 < x < \pi/2),
\end{align*}
appear in the fascinating identity  
\begin{align*} 
	\lim_{p \to \infty} \ehk{R_{p,d}} = 1 + m_d,
\end{align*}
established by Monsky and Gessel in \cite{gessel2010limit} building on ideas introduced in \cite{HanMonskySomeSH}. \\

We point out that both inequalties appearing in (\ref{WYineq}) are part of Conjecture \ref{WY Conjecture}.  The methods of \cite{HanMonskySomeSH} can be used to compute $\ehk{R_{p, d}}$ explicitly by hand --- though, the authors caution that these computations are highly involved and impractical in high dimension at present.  Using this approach, the second inequality in (\ref{WYineq}) has been established up to dimension $6$ by Yoshida \cite{yoshida2009small}.  Using different methods Trivedi has established the second inequality in all dimensions for large enough $p$ (see \cite{trivedi2021lower} and \cite{trivedi2023hilbert}).  \\

Conjecture \ref{WY Conjecture} is known to hold for rings with Hilbert-Samuel multiplicity at most $5$ by results in \cite{IanNewEst}, and is known in full generality for complete intersections by the work of Enescue and Shimomoto \cite{enescu2005upper}.  After introducing Conjecture \ref{WY Conjecture} in \cite{WY3Dim}, Watanabe and Yoshida proved the full conjecture up to dimension $4$.  This has since been extended to dimensions $5$ and $6$ by Aberbach and Enescue \cite{IanNewEst}. \\

The authors enthusiastically recommend the expository article \cite{huneke2014hilbert} for details about Hilbert-Kunz theory.  The recent article \cite{jeffries2023lower} contains additional results related to Conjecture \ref{WY Conjecture}, including a very interesting variant involving the $F$-signature.   \\

\subsection{Main Results}
\indt In this paper we establish that Conjecture \ref{WY Conjecture} is true up to dimension $7$: \\

\begin{thm} \label{WYDim7} \hfill \\
	
	The Watanabe-Yoshida Conjecture holds in dimension $d$, for all $d \le 7.$
\end{thm} 

\vspace{.25in}

Our approach is based on relating the Hilbert-Kunz multiplicity of $R$ with that of certain radical extensions, which we denote by $R_n.$\footnote{See section \ref{RadExtSection} for details about radical extensions.}  This approach to the problem is not new; radical extensions are employed in this context in \cite{aberbach2008lower} and the approach is further explored in \cite{IanNewEst} and \cite{CHDZ}.  Our contribution is to reformulate and combine the previous methods in a new way --- this leads to simplified arguments and improved bounds. \\

Our primary new tool is a function, $\Psifun{J}{R}{z^{t}},$ introduced and studied in section \ref{PhiDefConcave}, that interpolates between the Hilbert-Kunz multiplicities of the $R_n$'s.\footnote{This function was partly inspired by the convex function introduced by Blickle, Schwede, and Tucker in \cite{blickle2012f}.  Our function is related to theirs when $R$ is regular.}  The $\varphi$-function is increasing and equal to $\ehkideal{R}{J}$ when $t = 1$.  It enjoys a number of interesting analytic properties (see remark \ref{PhiPropertiesRemark}); in particular, we prove in Theorem \ref{Concave} that it is concave. \\

In section \ref{PhiBounds} we show how to combine the bounding methods of \cite{IanNewEst} and \cite{CHDZ} to produce non-trivial bounds for $\varphi$.  The resulting bound for $\ehk{R}$ appears in Corollary \ref{NoRootsBound}.  This lower bound for $\ehk{R}$ becomes more potent in Theorem \ref{BigBoundThm} where we involve multiple radical extensions.  
\begin{thm} (Theorem \ref{BigBoundThm}) \label{BigBoundThmIntro} \hfill \\

Suppose that $R$ is a complete, normal, local domain of positive characteristic and dimension $d$ with $\kappa = \ol{\kappa},$ and that $\unideal{x}$ is a minimal reduction of $\mf{m}_R$.  Let $\mu$ denote the number of minimal generators of $\mf{m}_R/\unideal{x}^{*}$ and suppose that $z = z_1, z_2 \dots, z_{\mu} \in \mf{m}_R \setminus \mf{m}_R^2$ map to a complete set of minimal generators of $\mf{m}_R/\unideal{x}^{*}.$ \\
\indt	Suppose that, for some $k < \mu,$ $z_1, z_2, \dots, z_{k+1}$ belong to a common minimal reduction of $\mf{m}_R$ in $R$ and for $i = 1, \dots, k,$ let 
	\begin{align*}
		S_i = R\big[z_2^{1/2}, z_3^{1/2}, \dots, z_{i+1}^{1/2}\big].	
	\end{align*}
\indt   Assume additionally that the radical extensions, $S_i,$ for $i = 1, \dots, k$, are all normal. Set $S = S_{k}.$ \\
\indt Then, for all $t \in [0, 1],$ and all $s \ge 0,$ \\
	\begin{align} \label{SinEqIntro}
		\ehk{S} \ge 1 - t + 2^k e\left(R\right) \left(\nu_{s} - \left(\mu - k - 1\right)\nu_{s-1} - k \nu_{s-1/2} - \nu_{s-t}\right)
	\end{align}
	and, hence
	\begin{align} \label{RinEqIntro}
		\ehk{R} \ge 1 - \frac{t}{2^k} +  e\left(R\right) \left(\nu_{s} - \left(\mu - k - 1\right)\nu_{s-1} - k \nu_{s-1/2} - \nu_{s-t}\right)
	\end{align}
\end{thm}

In the statement above, $z=z_1$ is used to generate the radical extensions that determine the $\varphi$-function, while the square roots of $z_2, \dots, z_k$ generate a sequence of quadratic extensions, $S_i.$  When $S=S_k$ is normal, we apply Corollary \ref{NoRootsBound} to get the inequality (\ref{SinEqIntro}) and then apply an improved version of Corollary 3.5 from \cite{CHDZ} (see remark \ref{PropPlus(v)}) to get the bound (\ref{RinEqIntro}).  

This bound represents a significant improvement, but it only works when we can adjoin enough of these square roots; in particular, $S=S_k$ needs to be normal in the setting above for our methods to work.  To address this, we prove a complementary bound in Lemma \ref{NotNormBnd}:
 \begin{lemma} (Lemma \ref{NotNormBnd}) \label{NotNormBndIntro} \hfill \\

Suppose $R$ and $z = z_1, z_2 \dots, z_{k+1} \in \mf{m}_R \setminus \mf{m}_R^2$ are as in the statement of Theorem \ref{BigBoundThmIntro}, and let 
	\begin{align*}
		S_i = R\big[z_2^{1/2}, z_3^{1/2}, \dots, z_{i+1}^{1/2}\big].	
	\end{align*}
Assume that $S_i$ is normal for $i=1, \dots, k-1,$ while $S_{k}$ is  not normal. \\
Then,
\begin{align*}
	\ehk{R} \ge 1 + \frac{1}{2^{k}}
\end{align*}
\end{lemma}

In section \ref{Dim7} we prove Conjecture \ref{WY Conjecture} in dimension $7$.  We provide a formula, computed using the methods of \cite{HanMonskySomeSH}, for $\ehk{R_{p, 7}}$\footnote{As far as the authors know, this is the first time this formula has appeared in print.}
\begin{align*}
	\ehk{R_{p, 7}} = \frac{332p^4+304p^2+192}{315p^4+273p^2+168}.
\end{align*} 

Using this, we establish the second inequality of Conjecture \ref{WY Conjecture} in dimension $7$.  To establish the remainder of the conjecture we utilize the bounds in Theorem \ref{BigBoundThm} with $k=1$ (justified by Lemma \ref{NotNormBnd}).  The critical bounds were obtained numerically, using \texttrademark{MATHEMATICA} --- they are listed in two tables in section \ref{NewDim7Bnds}. 
\newpage 
\section{Setup and Notation}
\subsection{Reducing to the Case of a Complete, Normal Domain with Algebraically Closed Residue Field} \label{reductionblurb}

\indt The techniques introduced in this paper are most naturally formulated for complete local normal domains with algebraically closed residue field.  Given our interest in Conjecture \ref{WY Conjecture}, where the local rings are only assumed to be formally unmixed, we begin by showing that we can make these simplifying assumptions.  
 
\begin{lemma} \label{reductionlemma} \hfill \\

	Suppose that $(R, \mf{m}_R, \kappa)$ is a formally unmixed Noetherian local ring of dimension $d \ge 2$ and characteristic $p>0.$  Assume that $R$ is not regular. \\
	
	 There exists a Noetherian local ring, $R'$ of characteristic $p$ and dimension $d$ which is a non-regular, complete normal domain with algebraically closed residue field such that 
	\begin{align*}
		\ehk{R'} \le \ehk{R}
	\end{align*}     
\end{lemma}
\begin{proof}
	Completing $R$ with respect to $\mf{m}_R$, does not change the lengths involved in the computation of Hilbert-Kunz multiplicity.  Thus we replace $R$ with it's completion --- and recall that, by assumption, the completion is unmixed. \\
	
Now choose an algebraic closure for the residue field, $\ol{\kappa}$ and consider the extension $R \xhookrightarrow{} R\widehat{\otimes}_{\kappa}\ol{\kappa},$ where $\widehat{\otimes}_{\kappa}$ denotes the completed tensor product over $\kappa$.\footnote{See section V.B.2 of \cite{serrelocal} for details about the completed tensor product --- and recall that $\ol{\kappa}$ is the directed colimit of finite extensions of $\kappa$.}  Noting that $\ol{\kappa}$ is faithfully flat over $\kappa$ and applying Cohen's Structure Theorem to $R$, we see that $R \widehat{\otimes}_{\kappa}\ol{\kappa}$ is the quotient of a Noetherian power series ring over $\ol{\kappa}$ and therefore itself complete, local and Noetherian.  We also see that the lengths involved in the Hilbert-Kunz function are unchanged if we replace $R$ with $R\widehat{\otimes}_{\kappa}\ol{\kappa}$.    Hence, we assume that $R$ is complete and unmixed with algebraically closed residue field.\footnote{This reduction, to the case where $R$ is complete with algebraically closed residue field, was first observed by Kunz in section 3 of his famous paper on the Frobenius endomorphism \cite{kunz1969characterizations}.} \\

If $R$ is not a domain then, by associativity, the Hilbert-Kunz multiplicity of $R$ is at least $2$.  In this case we set $R' = R_{p,d}$, where $R_{p,d}$ is the quadratic hypersurface of dimension $d$ and characteristic $p$ that features in the W-Y conjecture. Recall that $R_{p,d}$ is a non-regular complete normal domain satisfying $\ehk{R_{p, d}} < 2$ for $d\ge2$.\footnote{The inequality $\ehk{R_{p,d}}<2$ is easily verified from the identity $\ehk{R_{p,d}} + s(R_{p,d}) = 2$ by observing that $R_{p,d}$ is strongly $F$-Regular, and therefore $s(R_{p, d})>0$ --- where $s(R_{p, d})$ denotes the $F$-signature of $R_{p,d}.$} \\

We may now assume $R$ is a complete domain with algebraically closed residue field. Let $S$ denote the normalization of $R$ in it's quotient field.  $R$ is $F$-finite, and therefore it is excellent and the the extension $R \to S$ is finite.  Therefore $S$ is complete.  Moreover, since $R$ is Hensilian and $S$ is a domain, $S$ has a single maximal ideal, so it is local.  Letting $\mf{m}_S$ denote the maximal ideal of $S$, we have $$\ehk{S} = \ehkideal{S}{\mf{m}_S} \le \ehkideal{S}{\mf{m}_R S} = \ehk{R},$$ where the last equality follows from the fact that $S$ is the normalization of $R$, so it has rank $1$ as an $R$-module, by definition. \\

Now, it is possible that $S$ is regular, or equivalently that $\ehk{S} = 1$.  But, in this case $$\ehk{R} = \ehkideal{S}{\mf{m}_R S} \ge 2,$$ since $\mf{m}_R S \subset \mf{m}_S$ must be proper --- if $\mf{m}_R S = \mf{m}_S,$ then $\ehk{S} = \ehk{R}$ by the remarks above. As before, if $\ehk{R} \ge 2$, we can take $R' = R_{p,d}.$ \\

Otherwise, $S$ is a non-regular normal local domain with algebraically closed residue field, with $\dim S = \dim R,$ and $\ehk{S} \le \ehk{R}.$  Setting $R' = S$ we get the result as claimed.
\end{proof}
\subsection{Radical Extensions} \label{RadExtSection}
	\indt Throughout $(R, \mf{m}_R, \kappa)$ will denote a complete normal Noetherian local domain of positive characteristic $p > 0,$ with algebraically closed residue field, $\overline{\kappa} = \kappa.$  \\
	
	For a fixed minimal generator of $\mf{m}_R$, $z \in \mf{m}_R \setminus \mf{m}_R^2$ and $n \ge 1$ we let $R_n$ denote the {\it radical extension,} 
\begin{align*}
	R_n := \dfrac{R[Y]}{(Y^n - z)} = R[z^{1/n}],
\end{align*}
where $z^{1/n} \in R^{+}$ is a $n$th root of $z.$  \\

The reader can find additional information about radical extensions in Section 4 of \cite{aberbach2008lower} and Section 3 of \cite{CHDZ}.  For convenience, we use the notation $z^{a/b} = \left(z^{1/b}\right)^a.$
\begin{remark} \label{Rmk1} \hfill \\

These radical extensions are free: $R_n \cong R^{\oplus n}$ for each $n,$ with $R$-basis $1, z^{1/n}, z^{2/n}, \dots, z^{n-1/n}.$  Note that the $R_n$ also satisfy the following:
\begin{itemize}
	\item[(i)] $R[Y]$ is a domain, and $Y^n - z$ is clearly a regular element on $R[Y].$  In particular, if $R$ is CM or Gorenstein, then the same is true of all of the $R_n$.
	\item[(ii)] each $R_n$ is a local domain with maximal ideal $(\mf{m}_R, z^{1/n})R_n.$ \footnote{See Lemma \ref{PropPlus}.}
	\item[(iii)] For any $n, n', k \ge 1$ we have $z^{k/n} = z^{kn'/nn'}$ in $R_{nn'}.$
	\item[(iv)] $R/(z) \cong R_n/(z^{1/n})$ for all $n$. \\
	\item[(v)] For all $n, n' \ge 1,$ $$R_{nn'} \cong \dfrac{R_n [X]}{(X^{n'} - z^{1/n})}.$$ 
\end{itemize}
Notice that, if the $R_n$ are normal for all $n \ge 1$ (e.g. all the $R_n$'s are F-Regular), then all the statements above apply when $R$ is replaced by $R_n$ and $z$ is replaced by $z^{1/n}.$
\end{remark} \\  

\vspace{.5in}

\begin{remark} \label{CHDZrmk} \hfill \\

We are particularly interested in relating the Hilbert-Kunz multiplicity of $R$ to the Hilbert-Kunz multiplicities of these radical extensions.  This relationship is explored in \cite{IanNewEst} and \cite{CHDZ}, both of which are heavy influences on this paper.  Some known results in this direction are:   
	\begin{itemize}
		\item[(i)] We have $$\ehk{R} \le \ehk{R_n} \le \ehk{\ol{R}},$$ for all $n \ge 1,$ where $\ol{R} = R/(z).$ \\

Here, the first inequality follows because these radical extensions are flat, and the second inequality holds by filtration and the observation that $$R_n /(z^{1/n}) \cong \ol{R}$$ for all $n$. 
		\item[(ii)] When $z$ is part of a minimal reduction of $\mf{m}_R$, it is clear that the Hilbert-Samuel Multiplicities of $R$ and the $R_n$ are all equal, $$e(R) = e(R_n).$$
		\item[(iii)]  When $R$, and therefore all the $R_n$, are CM, the results of Polstra and Smirnov \cite{polstra_smirnov_2018} (see also \cite{CS20}, for a more general version of these results) shows that $$\displaystyle \lim_{n \to \infty} \ehk{R_n} = \ehk{\ol{R}}.$$		
	\end{itemize}
\end{remark}
\indt The following lemma consolidates several results for convenience, including a generalization of Proposition 3.4 of \cite{CHDZ} as well as some of it's consequences. \\

\begin{lemma} \label{PropPlus}
	Let $\left(R, \mf{m}_R, \kappa \right)$ be a complete normal Noetherian local domain of dimension $d$ with algebraically closed residue field.  Suppose that $\left(\ul{x}, z\right) = (x_1, \dots, x_{d-1}, z) \subset \mf{m}_R$ is a minimal reduction of $\mf{m}_R,$ where $z \in \mf{m}_R \setminus \mf{m}_R^2$ is a minimal generator of $\mf{m}_R$ and set $$S = R[z^{1/n}],$$ where $n > 1$ and $z^{1/n} \in R^{+}$ is an $n$th root of $z$. \\ 

	Then: \\
	\begin{itemize}
		\item[(i)] $S$ is a complete local domain with maximal ideal $\left(\mf{m}_R, z^{1/n} \right)S$. \\
		\item[(ii)] For each $i = 1, 2, \dots, n$ the ideal $\left(\ul{x}, z^{i/n}\right)S$ is a minimal reduction of $\left(\mf{m}_R, z^{i/n}\right)S$ in $S$. \\
		\item[(iii)] Therefore, the ideals $\left(\mf{m}_R, z^{i/n}\right)S$, for $i=1, 2, \dots, n$, are all integrally closed in $S$.
		\item[(iv)] When $R$ is $F$-finite of positive characteristic $p>0$, the Hilbert-Kunz multiplicities of the ideals $\left(\mf{m}_R, z^{i/n}\right)S$ satisfy the bound,
		\begin{align*}
			\ehkideal{S}{\left(\mf{m}_R, z^{i/n}\right) S} - \ehkideal{S}{\left(\mf{m}_R, z^{(i-1)/n}\right) S} \ge 1		
		\end{align*}
		for all $i=1, \dots, n.$ 
	\end{itemize}
\end{lemma}
\begin{proof} \hfill \\

	The fact that $\left(\mf{m}_R, z^{1/n}\right)S \subset S$ is maximal follows from the isomorphism $$S/\left(\mf{m}_R, z^{1/n}\right)S \cong R/\mf{m}_R \cong \kappa.$$  This is the only maximal ideal of $S$ because $S \subset R^{+}$ is a domain, $R \xhookrightarrow{} S \cong R[Y]/(Y^n-z)$ is finite and $R$ is Henselian.\footnote{Note that we are assuming $R$ is normal and $z \not\in \mf{m}_R^2$, and this ensures that $X^n - z$ is the minimal polynomial for $z^{1/n}$ over $R$ --- see remark 4.3 of \cite{aberbach2008lower} for details.}  \\

There is an isomorphism $S/\mf{m}_RS \cong \kappa[Y]/(Y^n),$ and therefore the $n$ ideals in the chain $$\mf{m}_R S = \left(\mf{m}_R, z^{n/n}\right)S \subset \left(\mf{m}_R, z^{n-1/n}\right)S \dots \subset \left(\mf{m}_R, z^{1/n}\right)S$$ are distinct and any ideal of $S$ that contains $\mf{m}_RS$ must equal $\left(\mf{m}_R, z^{i/n}\right)S$ for some $i=1, \dots, n.$ \\
	\indt Noting that $x_1, \dots, x_{d-1}, z^{1/n}$ form a system of parameters on $S,$ \cite{HSintegral} 11.2.9(2) gives
	\begin{align*}
		e\left(\left(\ul{x}, z^{i/n}\right)S, S\right) = e\left(\left(\ul{x}, \left(z^{1/n}\right)^i\right)S, S\right) = i e\left(\left(\ul{x}, z^{1/n}\right)S, S\right).
	\end{align*}

It follows that the multiplicities of the $n$ ideals in the chain $$\ol{\left(\ul{x}, z^{n/n}\right)S} \subset \ol{\left(\ul{x}, z^{n-1/n}\right)S} \dots \subset \ol{\left(\ul{x}, z^{1/n}\right)S}$$ are all different, and hence, by \cite{HSintegral} proposition 11.2.1, they are all distinct integrally closed ideals of $S$.  \\

Since these ideals all contain $\mf{m}_RS$, 
	\begin{align*}
		\mf{m}_RS = \ol{\left(\ul{x}, z\right)}S \subset \ol{\left(\ul{x}, z^{n/n}\right)S} \subset \ol{\left(\ul{x}, z^{i/n}\right)S},		
	\end{align*}	
	the observations above show that we must have $\ol{\left(\ul{x}, z^{i/n}\right)S} = \left(\mf{m}_R, z^{i/n}\right)S$ for each $i = 1, \dots, n,$.  This establishes (ii) and (iii) follows immediately. \\
	\vspace{.2in}

For each $i = 1, \dots, n$ and all $q = p^e,$ we have $$\len{S}{\frac{\left(\left(\mf{m}_R, z^{i/n}\right)S\right)^{\fbp{q}}}{\left(\left(\mf{m}_R, z^{(i-1)/n}\right)S\right)^{\fbp{q}}}} =  \len{S}{\frac{\left(\mf{m}_R, z^{i/n}\right)S^{1/q}}{\left(\mf{m}_R, z^{(i-1)/n}\right)S^{1/q}}} \ge \rank{S^{1/q}} = q^{d}$$ by corollary 2.2 of \cite{CHDZ}.\footnote{Note that $S$ is complete with perfect residue field, so $S$ is $F$-finite.}  Dividing and taking the limit $q \to \infty$ gives the inequality.
\end{proof} 

\begin{remark} \label{PropPlus(v)}
	In (iv) above we have 
	\begin{align*}
		\ehkideal{S}{\mf{m}_R S} = n \ehk{R},
	\end{align*}
	by Theorem 2.7 of \cite{watanabe2000hilbert}.  Applying the inequality in Lemma \ref{PropPlus} (iv) to the chain of ideals, $$\mf{m}_R S = \left(\mf{m}_R, z^{n/n}\right)S \subset \left(\mf{m}_R, z^{n-1/n}\right)S \dots \subset \left(\mf{m}_R, z^{1/n}\right)S$$ yields the useful bound,
	\begin{align*}
		\ehk{R} - 1 \ge \frac{\ehk{S} - 1}{n}.
	\end{align*}
	This is the inequality of Corollary 3.5 in \cite{CHDZ} generalized to our setting. 
\end{remark}
\vspace{.1in}

\newpage

\section{The Convex $\varphi_J$ Function} \label{PhiSection}
\indt Continuing to let $(R, \mf{m}_R, \kappa)$ denote a complete, normal local domain with $\kappa = \ol{\kappa},$ and taking $z \in \mf{m}_R \setminus \mf{m}_R^2$ as above, suppose that $J \subset R$ is an $\mf{m}_R$-primary ideal.  We are going to introduce a function, $\Psifun{J}{R}{z^t},$ defined for real valued $t$, that interpolates the Hilbert-Kunz multiplicities of ideals in the radical extensions $R_n$ discussed above.  We will see that this function is very well behaved and that it fits conveniently in the framework of \cite{IanNewEst} and \cite{CHDZ}.  In particular, the results of \cite{IanNewEst} naturally generalize to give explicit numerical bounds on the values of the $\varphi$ function.  In combination with ideas introduced in \cite{CHDZ} this can be used to obtain non-trivial bounds for the Hilbert-Kunz multiplicities of ideals in $R$.

\vspace{.5in}

\subsection{Definition and Concavity} \label{PhiDefConcave}
\begin{lemma} \hfill \\

Suppose that $a,b, c, d \ge 1$ are integers with $a/b = c/d.$  Then
\begin{align*}
	\dfrac{1}{b} \, \, \ehkideal{R_{b}}{\left(J, z^{a/b}\right)R_{b}} = \dfrac{1}{d} \, \, \ehkideal{R_{d}}{\left(J, z^{c/d}\right)R_{d}}. 
\end{align*} 
\end{lemma}  
\begin{proof} \hfill \\

By Theorem 2.7 of \cite{watanabe2000hilbert} we have:
\begin{small}	
	\begin{align*}
	\ehkideal{R_{bd}}{\left(J, z^{ad/bd}\right)R_{bd}} &= \ehkideal{R_{bd}}{\left(J, z^{a/b}\right)R_{bd}} \\
	                                     &= [QF(R_{bd}) \, : \, QF(R_b)]\ehkideal{R_{b}}{\left(J, z^{a/b}\right)R_{b}} \\
	                                     &= \dfrac{[QF(R_{bd}) \, : \, QF(R)]}{[QF(R_{b}) \, : \, QF(R)]}\ehkideal{R_{b}}{\left(J, z^{a/b}\right)R_{b}} \\
	                                     &= \dfrac{bd}{b} \,\, \ehkideal{R_{b}}{\left(J, z^{a/b}\right)R_{b}} \\
	                                     &= d \,\, \ehkideal{R_{b}}{\left(J, z^{a/b}\right)R_{b}}
	\end{align*}
\end{small}
	Similarly,
	\begin{small}
	\begin{align*}
		\ehkideal{R_{bd}}{\left(J, z^{bc/bd}\right)R_{bd}} = b \,\, \ehkideal{R_{d}}{\left(J, z^{c/d}\right)R_{d}}
	\end{align*}
	\end{small}
	and by assumption $ad = bc,$ so we have 
	\begin{small}
	\begin{align*}
		d \,\, \ehkideal{R_{b}}{\left(J, z^{a/b}\right)R_{b}} = b \,\, \ehkideal{R_{d}}{\left(J, z^{c/d}\right)R_{d}}.
	\end{align*}
	\end{small}
\end{proof}
\indt Therefore, the following function is well defined: \\

\begin{defn} \hfill \\

For $a, b > 0,$
\begin{align*}
	\Psifun{J}{R}{z^{a/b}} := \dfrac{1}{b} \, \, \ehkideal{R_{b}}{\exideal{J}{z^{a/b}}{R_d}}
\end{align*}
We also set $\Psifun{J}{R}{z^{0}} = 0.$
\end{defn}
\begin{remark} \hfill \\

It is immediate from this definition that $\Psifun{J}{R}{z^{a/b}}$ is an increasing function of $a/b$ and that $\PsifunZ{1} = \ehkideal{R}{\left(J, z\right)}.$ \\
\end{remark}

\begin{remark} \hfill \\

Many of the nice analytic properties of $\varphi_J$ come down to the fact that it is a convex function --- which is proven in lemma \ref{Concave} below.  The convexity of $\varphi_J$ is essentially equivalent to the existence of the surjections appearing at the start of the proof of lemma \ref{concave1}.   
\end{remark} \\

The next result establishes several useful properties of $\varphi_J$ which, among other things, allow us to extend this definition to $\varphi_{J}\left(R; \, z^t\right)$ for real values of $t\ge 0.$ \\
\begin{lemma} \label{concave1}  \hfill \\ 

For $a, b > 0,$
\begin{small}
\begin{align*}
	0 \le \PsifunZ{(a+1)/b} - \PsifunZ{a/b} \le \PsifunZ{a/b} - \PsifunZ{(a-1)/b} \le \PsifunZ{1/b} \le \dfrac{1}{b} \,\, \ehkideal{\ol{R}}{J\ol{R}}
\end{align*}
\end{small}
where $\ol{R} = R/(z)$ and $\ehkideal{\ol{R}}{J\ol{R}}$ denotes the Hilbert-Kunz Multiplicity of the ideal $J\ol{R}$ computed in $\ol{R}.$ \\

In particular, the quantity
\begin{small}
\begin{align*} 
	\dfrac{1}{\ell} \left[\PsifunZ{(a+\ell)/b} - \PsifunZ{a/b}\right]
\end{align*}
\end{small}
is decreasing as a function of both $a \ge 0$ and $\ell \ge 1$.
\end{lemma}
\begin{proof} \hfill \\

To see the first two inequalities, note that for any $q = p^e$ there are surjections
	\begin{small}
	\begin{align*}
		\dfrac{\left(J^{\fbp{q}}, z^{(a-1)q/b} \right)R_b}{\left(J^{\fbp{q}}, z^{aq/b} \right)R_b} \to \dfrac{\left(J^{\fbp{q}}, z^{aq/b} \right)R_b}{\left(J^{\fbp{q}}, z^{(a+1)q/b} \right)R_b} \to 0 
	\end{align*}
	\end{small}
	and 
	\begin{small}
\begin{align*}
		\dfrac{R_b}{\left(J^{\fbp{q}}, z^{q/b} \right)R_b} \to \dfrac{\left(J^{\fbp{q}}, z^{(a-1)q/b} \right)R_b}{\left(J^{\fbp{q}}, z^{aq/b} \right)R_b} \to 0 
	\end{align*}
	\end{small}

The last inequality follows from the fact that $R_b/(z^{1/b}) \cong R/(z)$ and the standard filtration argument. \\

For the final statement, that the indicated combination is decreasing in $a$ is immediate from the inequalities above.  To see that it also decreases in $\ell$, suppose that $1 \le \ell.$ \\

By the inequalities established above,
	\begin{small}
	\begin{align*}
		\PsifunZ{(a+\ell+1)/b}- \PsifunZ{(a+\ell)/b} \le \PsifunZ{(a+h)/b} - \PsifunZ{(a+h-1)/b}
	\end{align*}
	\end{small}
	for every $h$ with $1 \le h \le \ell.$  Therefore,
	\begin{small}
	\begin{align*}
		\ell &\left[\PsifunZ{(a+\ell+1)/b} -\PsifunZ{a/b}\right] \\
		&= \ell \left[\PsifunZ{(a+\ell+1)/b}- \PsifunZ{(a+\ell)/b}\right] + \ell \left[\PsifunZ{(a+\ell)/b}- \PsifunZ{a/b}\right] \\
		&\le \sum_{h=1}^{\ell}\left(\PsifunZ{(a+h)/b} - \PsifunZ{(a+h-1)/b}\right) + \ell \left[\PsifunZ{(a+\ell)/b}- \PsifunZ{a/b}\right] \\
		&= (\ell + 1 )\left[\PsifunZ{(a+\ell)/b} - \PsifunZ{a/b}\right]
	\end{align*}
	\end{small}
\end{proof}

In what follows, we will continue to use $\ol{R}$ to denote the quotient $R/(z).$
\begin{lemma} \label{lip}  \hfill \\

For $a, c \ge 0$ and $b, d > 0,$ 
\begin{small}
\begin{align*}
	\left| \PsifunZ{a/b} - \PsifunZ{c/d} \right| \le  \left|\dfrac{a}{b} - \dfrac{c}{d} \right| \ehkideal{\ol{R}}{J\ol{R}}
\end{align*}
\end{small}
\end{lemma}  
\begin{proof} \hfill \\

First note that the two sides are both zero if $a/b = c/d$.  Otherwise we may choose a common denominator and relabel so that the numbers under consideration are $a/b$ and $c/b$, with $a/b > c/b.$ \\

We have, 
\begin{small}
\begin{align*}
	\PsifunZ{a/b} - \PsifunZ{c/b} &= \sum_{\ell = 1}^{a-c} \left[\PsifunZ{(c+\ell)/b} - \PsifunZ{(c+\ell-1)/b}\right] \\
	           &\le \sum_{\ell = 1}^{a-c} \frac{1}{b}  \ehkideal{\ol{R}}{J\ol{R}}  
\end{align*}
\end{small}
by lemma \ref{concave1}. \\

This last expression is equal to $\dfrac{a-c}{b} \ehkideal{\ol{R}}{J\ol{R}}$, so the inequality is established.
\end{proof}

We are now in position to extend our definition of $\varphi_J$ to the positive real numbers. \\
\begin{defn} \hfill \\

For any $t \in \mathbb{R}_{\ge 0},$ let $r_k \to t$ be any sequence of positive rational numbers converging to $t,$ and define 
\begin{align*}
	\PsifunZ{t} := \displaystyle \lim_{k\to\infty} \PsifunZ{r_k}
\end{align*}
This limit exists and is independent of the choice of $r_k \to t$ by lemma \ref{lip}.
\end{defn}

\begin{lemma} \label{Concave} \hfill \\

For any $t_1, t_2 \in \mathbb{R}_{\ge 0},$ and any $\alpha \in [0, 1],$ 
\begin{align*}
	\PsifunZ{\alpha t_1 + (1-\alpha)t_2} \ge \alpha \PsifunZ{t_1} + (1-\alpha) \PsifunZ{t_2}
\end{align*}
That is, $\varphi_J,$ is upper convex (i.e. concave down).
\end{lemma}
\begin{proof} \hfill \\

It is enough to prove this when $t_1, t_2$ and $\alpha$ are all rational. Also, if $t_1 = t_2$ or if $\alpha = 0,$ then the statement is trivial, so we will assume $t_1 > t_2$ and $\alpha > 0.$ \\

Write $\alpha = a/b$ and $t_1 = (c + w)/d$ and $t_2 = c/d$ --- with $c \ge 0$, $a, b, d, w \ge 1$ and $b \ge a$, since $\alpha \in [0, 1]$.\\ 
Note that, with these assignments,
\begin{small} 
\begin{align*}
	\alpha t_1 + (1-\alpha) t_2 = \dfrac{a(c+w)}{bd} + \dfrac{(b-a)c}{bd} = \dfrac{aw + bc}{bd}
\end{align*}
\end{small}

The final statement of lemma \ref{concave1} gives
\begin{small}
\begin{align*}
	\dfrac{1}{aw} \left[ \PsifunZ{(aw + bc)/bd} - \PsifunZ{bc/bd}\right] \ge  \dfrac{1}{bw} \left[ \PsifunZ{(bw + bc)/bd} - \PsifunZ{bc/bd}\right]
\end{align*}
\end{small}
Simplifying this gives,
\begin{small}
\begin{align*}
	\PsifunZ{(aw + bc)/bd} - \PsifunZ{c/d} \ge  \dfrac{a}{b} \left[ \PsifunZ{(w + c)/d} - \PsifunZ{c/d}\right]
\end{align*}
\end{small}
Adding $\PsifunZ{c/d}$ to both sides, we get
\begin{small}
\begin{align*}
	\PsifunZ{(aw + bc)/bd} \ge \dfrac{a}{b}\PsifunZ{(w+c)/d} + \dfrac{b-a}{b}\PsifunZ{c/d}
\end{align*}
\end{small}
as desired.
\end{proof} 
\begin{remark} \label{PhiPropertiesRemark} \hfill \\

Several immediate consequences of the results above are worth mentioning. \\
	\begin{itemize}
		\item[(i)] Lemma \ref{lip} says that $\PsifunZ{t}$ is Lipschitz continuous with constant $\ehkideal{\ol{R}}{J\ol{R}}$. \\
		\item[(ii)] Convex functions are differentiable almost everywhere and absolutely continuous.  In particular, the left and right derivatives of $\PsifunZ{t}$, $$D_{-}\PsifunZ{t} = \displaystyle \lim_{s \to 0^{+}} \dfrac{\PsifunZ{t} - \PsifunZ{t-s}}{s}$$ and $$D_{+}\PsifunZ{t} = \displaystyle \lim_{s \to 0^{+}} \dfrac{\PsifunZ{t+s} - \PsifunZ{t}}{s},$$ exist for all $t>0,$ and are equal almost everywhere.  \\
		\item[(iii)] Moreover, $D_{-}\PsifunZ{t}$ and $D_{+}\PsifunZ{t}$ are continuous from the left and right, respectively, and are thus completely determined by their values at rational numbers (see, e.g., Theorem 2.6 of \cite{hug2020lectures}).  
		\item[(iv)] It follows from lemma \ref{Concave} that both $D_{-}\PsifunZ{t}$ and $D_{+}\PsifunZ{t}$ are decreasing functions of $t$ and, when $\ol{R}$ is generically reduced and formally  equidimensional, the results of \cite{CS20} show that $$D_{+}\PsifunZ{0} = \ehkideal{\ol{R}}{J\ol{R}}.$$ \\
		\item[(v)] The derivative $\varphi_{J}' \left(R; t\right)$ is defined a.e., and for all $t > t_0 \ge 0,$ there is an equality
		\begin{align*}
			\PsifunZ{t} = \PsifunZ{t_0} + \displaystyle \int_{t_0}^t \varphi_{J}' \left(R; \tau\right) d \tau		
		\end{align*}
		\item[(vi)]  It is also not difficult to show that $\varphi_{J}' \left(R; z^t\right)$ is of bounded variation, and that the second derivative $\varphi_{J}'' \left(R; z^t\right)$ exists a.e. \\
		\indt  Note though, that there appears to be no reason to expect that $\varphi_{J}' \left(R; z^t\right)$ is continuous.  For example, while the discussion above shows that the integral $\int \varphi_{J}'' \left(R; z^t\right) dt$ is well defined, the relationship between $\varphi_{J}' \left(R; z^t\right)$ and $\int \varphi_{J}'' \left(R; z^t\right) dt$ is not at all clear. 
	\end{itemize}
	 
\end{remark}
\subsection{Bounds Involving $\PsifunZ{t}$} \label{PhiBounds}
\indt The techniques of \cite{IanNewEst} naturally extend to give explicit bounds for $\PsifunZ{t}.$ \\

\begin{prop} \label{PropShop} \hfill \\

Suppose that $R$ is a complete normal, local domain of dimension $d$  with algebraically closed residue field, and continue to let $J \subset R$ be an $\mf{m}_R$-primary ideal.  Let $z \in \mf{m}_R \setminus \mf{m}_R^2$ be a fixed minimal generator of the maximal ideal.  Let $\unideal{x} \subset R$ denote an $\mf{m}_R$-primary parameter ideal such that $\unideal{x}^{*} \subset \left(J, z\right)$ and such that $z \in \overline{\unideal{x}}.$  Suppose that $z = z_1, z_2 \dots, z_{r} \in \left(J, z\right)$ generate $\left(J, z\right)/\unideal{x}^{*}$.  \\
\vspace{.1in}

For each $z_i, $ let $t_i = f_{\unideal{x}}(z_i),$ where
	\begin{align*}
		f_{\unideal{x}}(y) = \lim_{k \to \infty} \dfrac{\sup \left(\ell \, \, \, \, | \, \, \, \,  y^k \in \unideal{x}^\ell \right)}{k}	
	\end{align*} 
	and for real valued $s,$ let \footnote{Observe that $\nu_s$ is a non-negative {\it increasing} function with $\nu_s = 1$ for $s \ge d$ and $\nu_s = 0$ for all $s \le 0.$  Several interesting formulas are known for $\nu_s$ --- e.g. see \cite{cho2020volume, marichal2006slices} --- including the following, which hold for $s \ge 0$: 
	\begin{align*}
		\nu_s &= \sum_{j=0}^{\floor{s}} \left(-1 \right)^j \dfrac{(s-j)^d}{j! ( d-j)!} = \frac{2}{\pi} \int_{0}^s\int_{0}^{\infty} \left(\frac{\sin u}{u}\right)^d \cos \left((d-2w)u\right) du \,\, dw
	\end{align*}}
	\begin{align*}
		\nu_{s} &:= \mathrm{vol} \{ (x_1, \dots, x_d) \in [0, 1]^d \, \, | \, \, \sum_{i=1}^d x_i \le s \}. \\
	\end{align*}
	
Then, for all $t \in [0, 1],$ and all $s \ge 0,$
	\begin{align*}
		\PsifunZ{t} \ge e\left(R, \unideal{x} \right) \left( \nu_s - \sum_{k = 2}^{r} \nu_{s - t_k}  - \nu_{s - t} \right).		
	\end{align*}	
\end{prop}
\begin{proof} \hfill \\

For rational $k/n \in [0, 1],$ the proof of Theorem 3.2 in \cite{IanNewEst} applies to $\left(\unideal{x}^{*}\right)R_n \subset \left(J, z^{k/n} \right)R_n,$ and gives the bound \footnote{Note that the conditions of Lemma \ref{PropPlus} are satisfied and, therefore, $R$ and all of the $R_n$'s are complete domains, and, in particular, unmixed.  The statement of Theorem 3.2 in \cite{IanNewEst} assumes that $\left(\left(\unideal{x}^*\right)R_n\right)^* = \left(\unideal{x}R_n\right)^* \subset \left(J, z^{k/n} \right)R_n,$ but the argument works if we only assume the weaker inclusion $\left(\unideal{x}^* \right) R_n \subset \left(J, z^{k/n} \right)R_n.$  Also, note that $z^{k/n}, z_2, \dots, z_r$ are clearly generators for $(J, z^{k/n})R_n / \left(\unideal{x}^{*}\right)R_n$ and that the proof in \cite{IanNewEst} goes through even if they are not minimal generators.  With these (possibly) weaker assumptions, the strength of Theorem 3.2 in \cite{IanNewEst} is slightly weakened as well. }
\begin{align*}
	\ehkideal{R_{n}}{(J, z^{k/n})R_n} \ge e\left(R_n, \unideal{x} \right) \left(\nu_{s} -\sum_{k = 1}^{r} \nu_{s-t_i} \right),
\end{align*}
with\footnote{This holds because $R \xhookrightarrow{} R_n$ is split and therefore $z_i^k \in \unideal{x}^\ell R_n \Longleftrightarrow z_i^k \in \unideal{x}^\ell$ in $R.$}
 $$t_i = f_{\unideal{x}}(z_i) = f_{\unideal{x}R_n} (z_i) ,$$  for $i = 2, \dots, \mu,$ and (see \cite{HSintegral} for details) $$t_1 = f_{\unideal{x}R_n} \left(z^{k/n}\right) \ge k/n.$$  
Putting this into the inequality above, and using the fact that $e\left(\unideal{x}R_n\right)
= n e(R, \unideal{x}),$ we get \footnote{Recall that $\nu_{s}$ is an {\it increasing} function of $s$.}
\begin{align*}
	\ehkideal{R_{n}}{(J, z^{k/n})R_n} &\ge n \, e\left(R, \unideal{x} \right) \left(\nu_{s} -\sum_{k = 2}^{r} \nu_{s-t_i} - \nu_{s-t_1} \right) \\
		&\ge n \, e\left(R, \unideal{x} \right) \left(\nu_{s} -\sum_{k = 2}^{r} \nu_{s-t_i} - \nu_{s-k/n} \right)	
\end{align*}
Therefore, for any $t = k/n \in [0,1],$ we have
\begin{align*}
	\PsifunZ{t} \ge e\left(R, \unideal{x} \right) \left(\nu_{s} -\sum_{k = 2}^{r} \nu_{s-t_i} - \nu_{s-t} \right)	
\end{align*}
and the claim holds for arbitrary $t \in [0, 1]$ by continuity.
\end{proof}

Now we specialize to the case where $J = \mf{m}_R$, so that we can use the ideas introduced in \cite{CHDZ}. \\

Since this setup will be our primary focus going forward, we introduce the following notation,
\begin{align*}
	\varphi(t) := \Psifun{\mf{m}_R}{R}{z^t}
\end{align*}
to make things easier to read. \\

\begin{lemma} \label{CHDZIne} \hfill \\

We have, for all $t \ge 1,$
\begin{align*}
	\varphi(t) = \varphi(1) = \ehk{R}.
\end{align*}

Assume that the minimal generator, $z \in \mf{m}_R \setminus \mf{m}_R^2$, belongs to a minimal reduction of $\mf{m}_R.$ \footnote{This is well known to be a generic condition --- see Theorem 8.6.6 of \cite{HSintegral}.} \\

Then, for all $0 \le t_0 \le t \le 1,$ 
	\begin{align*}
		\varphi(t) \ge \varphi(t_0) + t - t_0
	\end{align*} 
\end{lemma} 
\begin{proof} \hfill \\

The fact that $\varphi(1) = \ehk{R}$ is immediate from the definition.  The other equality, $$\varphi(t) = \varphi(1)$$ for $t\ge 1,$ follows from the fact that $\left(\mf{m}_R, z^{k/n}\right)R_n = \mf{m}_R R_n,$ for any $k/n \ge 1,$ along with the observation that $\ehkideal{R_n}{\mf{m}_R R_n} = n \ehk{R}.$ \\

The second claim of the lemma is a consequence of lemma \ref{PropPlus}: by \ref{PropPlus} (iii) the ideals $\left(\mf{m}_R, z^{k/n}\right)R_n,$ with $1 \le k \le n,$ are all integrally closed; and \ref{PropPlus} (iv), shows that their relative Hilbert-Kunz multiplicities are at least $1$, $$\ehkideal{R_n}{\left(\mf{m}_R, z^{k/n}\right) R_n} - \ehkideal{R_n}{\left(\mf{m}_R, z^{(k-1)/n}\right) R_n} \ge 1.$$ 

Assume for the moment that $t$ and $t_0$ are rational and write $t_0 = k/n$ and $t = (k+\ell)/n,$ with $0 \le k \le k +\ell \le n.$  Applying the inequality above repeatedly, we have
\begin{small}
\begin{align*}
	\ehkideal{R_n}{\left(\mf{m}_R, z^{(k+\ell)/n}\right) R_n} &- \ehkideal{R_n}{\left(\mf{m}_R, z^{k/n}\right) R_n} \\ 
	&= \sum_{j=1}^{\ell}  \,\, \ehkideal{R_n}{\left(\mf{m}_R, z^{(k+j)/n}\right) R_n} - \ehkideal{R_n}{\left(\mf{m}_R, z^{(k+j-1)/n}\right) R_n} \\
	&\ge \ell.
\end{align*}
\end{small}
Dividing by $n$ and rearranging things a bit gives
\begin{align*}
	\frac{1}{n}\ehkideal{R_n}{\left(\mf{m}_R, z^{(k+\ell)/n}\right) R_n} \ge \frac{1}{n}\ehkideal{R_n}{\left(\mf{m}_R, z^{k/n}\right) R_n} + \frac{\ell}{n}.
\end{align*}
Writing this in terms of $\varphi, t,$ and $t_0,$ we get the desired inequality,
\begin{align*}
	\varphi(t) \ge \varphi(t_0) + t-t_0.
\end{align*}
This extends to general $t$ and $t_0$ by continuity.
\end{proof}
\begin{remark} \hfill \\

The proof of lemma \ref{CHDZIne} shows that $\varphi'(t) \ge 1,$ a.e. for $t \in [0, 1],$ and that $\varphi'(t) = 0$ for $t>1.$ \\
\end{remark} 

\vspace{.25in}

Putting Lemma \ref{CHDZIne} together with the bound in Proposition \ref{PropShop} we get a the following bound for $\varphi$. \\

\begin{cor} \label{NoRootsBound} \hfill \\

Assume $R$ is a complete, normal local domain of characteristic $p>0$ and dimension $d$ with algebraically closed residue field and that $z \in \mf{m}_R/\mf{m}_R^2$ belongs to a minimal reduction of $\mf{m}_R.$  Suppose that $\unideal{x}$ is an $\mf{m}_R$-primary parameter ideal with $z \in \ol{\unideal{x}}$, and assume that $z = z_1, z_2 \dots, z_{r} \in \mf{m}_R,$ generate $\mf{m}_R/\unideal{x}^{*}.$  As in Proposition \ref{PropShop}, let $t_i =  f_{\unideal{x}}(z_i)$ for $i = 2, 3, \dots, r.$ \\

Then for all $0 \le t_0 \le t \le 1,$ and all $s \ge 0,$ 
	\begin{align*}
		\varphi(t) \ge  t - t_0 + e\left(R, \unideal{x} \right) \left(\nu_{s} - \sum_{k = 2}^{r} \nu_{s - t_k} - \nu_{s-t_0}\right).
	\end{align*}

In particular, 
	\begin{align} \label{NoRootsInEq}
			\ehk{R} = \varphi(1) \ge  1 - t + e\left(R, \unideal{x} \right) \left(\nu_{s} -  \sum_{k = 2}^{r} \nu_{s - t_k} - \nu_{s-t}\right) 
	\end{align}
	for all $s \ge 0,$ and all $t \in [0, 1].$
\end{cor} 

\vspace{.25in}

In order to get the most out of these inequalities, we will need to choose the $z$'s optimally, so that we can adjoin some additional roots.

\begin{thm} \label{BigBoundThm} \hfill \\

Suppose that $R$ is a complete, normal, local domain of dimension $d$ with $\kappa = \ol{\kappa},$ and that $\unideal{x}$ is a minimal reduction of $\mf{m}_R$.  Let $\mu$ denote the number of minimal generators of $\mf{m}_R/\unideal{x}^{*}$ and suppose that $z = z_1, z_2 \dots, z_{\mu} \in \mf{m}_R \setminus \mf{m}_R^2$ map to a complete set of minimal generators of $\mf{m}_R/\unideal{x}^{*}.$ \\

Suppose that, for some $k < \mu,$ $z_1, z_2, \dots, z_{k+1}$ belong to a common minimal reduction of $\mf{m}_R$ in $R$ and for $i = 1, \dots, k,$ let \footnote{When treating the dimension $7$ case in the next section, we use this result with $k=1$ --- this turns out to give a more streamlined argument than taking larger $k$.}
	\begin{align*}
		S_i = R\big[z_2^{1/2}, z_3^{1/2}, \dots, z_{i+1}^{1/2}\big].	
	\end{align*}

Assume additionally that the radical extensions, $S_i,$ for $i = 1, \dots, k$, are all normal. Set $S = S_{k}.$ \\

Then, for all $t \in [0, 1],$ and all $s \ge 0,$ \\
	\begin{align} \label{SinEq}
		\ehk{S} \ge 1 - t + 2^k e\left(R\right) \left(\nu_{s} - \left(\mu - k - 1\right)\nu_{s-1} - k \nu_{s-1/2} - \nu_{s-t}\right)
	\end{align}
	and, hence
	\begin{align} \label{RinEq}
		\ehk{R} \ge 1 - \frac{t}{2^k} +  e\left(R\right) \left(\nu_{s} - \left(\mu - k - 1\right)\nu_{s-1} - k \nu_{s-1/2} - \nu_{s-t}\right).
	\end{align}
\end{thm}
\begin{proof} \hfill \\

Inequality (\ref{SinEq}) is simply inequality (\ref{NoRootsInEq}) of Corollary \ref{NoRootsBound} applied inside $S$: By lemma \ref{PropPlus} (i), $S$ is a complete, local, normal domain with maximal ideal $$\mf{m}_S = \left(\mf{m}_R, z_2^{1/2}, \dots, z_{k+1}^{1/2}\right)S,$$ and $z=z_1, z_2^{1/2}, \dots, z_{k+1}^{1/2}$ clearly generate $\mf{m}_S/\left(\unideal{x}^{*}\right)S,$ and therefore they map to generators in $\mf{m}_S/\left(\unideal{x}S\right)^*$, since $\left(\unideal{x}^{*}\right)S \subset \left(\unideal{x}S\right)^{*}.$ \footnote{Also, recall that $z=z_1 \in \mf{m}_S \setminus \mf{m}_S^2$ is part of a minimal reduction of $\mf{m}_S$ in $S$, by construction.} \\   

Noting that $S$ is free of rank $2^k$ over $R$ we conclude that $$e\left(S, \unideal{x}S\right) = 2^ke(R).$$  

Furthermore, $$f_{\unideal{x}S}\left(z_i^{1/2}\right) \ge 1/2$$ for $i = 2, \dots, k+1,$ since $\left(z_i^{1/2}\right)^2 = z_i \in \ol{\unideal{x}}S \subset \ol{\unideal{x}S}$ (\cite{HSintegral}), while for $i \ge k+2,$ it is clear that $$f_{\unideal{x}S}(z_i) \ge 1.$$

Applying \ref{PropPlus} (v) for each of the $k$-square roots generating $S$ over $R$, we get $$\ehk{R} - 1 \ge \dfrac{\ehk{S} - 1}{2^k}.$$  Combined with (\ref{SinEq}), this gives inequality (\ref{RinEq}).
\end{proof}

\newpage

\begin{remark} \hfill \\

The methods above can be modified to get a bound analogous to (\ref{RinEq}) that applies more generally to $\PsifunZ{t}$ --- at least for certain $J.$  Here we are concerned primarily with bounds on $\ehk{R},$ so we have saved these considerations for a later project. 
\end{remark} \\

Before we move on to the specific problem of bounding Hilbert-Kunz multiplicity in dimension $d = 7,$ we need to discuss what happens when we adjoin square roots in Theorem \ref{BigBoundThm} but the resulting ring, $S = S_k,$ is no longer normal. \\

As above, assume $R$ is a complete, normal local domain of positive characteristic $p>0$ and dimension $d \ge 2$ with algebraically closed residue field and that $z \in \mf{m}_R/\mf{m}_R^2$ belongs to a minimal reduction of $\mf{m}_R.$  Letting $S = R[z^{1/2}],$ we have seen that $S$ is automatically a complete local domain. \\

Suppose though, that $S$ is not normal.  Now, $S$ is a free $R$-module of rank $2$, so it is automatically $S_2.$  Thus, by Serre's criterion there is a prime $\mf{p} \in \spec{S}$ of height $1$ such that $S_{\mf{p}}$ is not regular.  $S$ is complete and equidimensional so $$\dim S/\mf{p} + \height{\mf{p}} = \dim S,$$ and therefore, $$\ehk{S} \ge \ehk{S_{\mf{p}}}$$ by Corollary 3.8 of \cite{kunz1976noetherian}.  It is known that singular local rings of dimension $1$ have Hilbert-Kunz multiplicity of at least $2$, and therefore, $$\ehk{S} \ge 2.$$  Using the inequality in remark \ref{PropPlus(v)}, we conclude that $$\ehk{R} \ge 1 + \dfrac{\ehk{S} - 1}{2} = 1 + \frac{1}{2}.$$ \\
This brings us to the following very useful result. \\

\begin{lemma} \label{NotNormBnd} 

Suppose $R$ and $z = z_1, z_2 \dots, z_{k+1} \in \mf{m}_R \setminus \mf{m}_R^2$ are as in the statement of Theorem \ref{BigBoundThm}, and let 
	\begin{align*}
		S_i = R\big[z_2^{1/2}, z_3^{1/2}, \dots, z_{i+1}^{1/2}\big].	
	\end{align*}
Assume that $S_i$ is normal for $i=1, \dots, k-1,$ while $S_{k}$ is  not normal. \\
Then,
\begin{align*}
	\ehk{R} \ge 1 + \frac{1}{2^{k}}
\end{align*}
\end{lemma}
\begin{proof}

The argument above shows that $\ehk{S_{k-1}} \ge 1 + \frac{1}{2},$ and, as in the proof of Theorem \ref{BigBoundThm}, repeated application of the inequality in remark \ref{PropPlus(v)} gives $$\ehk{R} - 1 \ge \frac{\ehk{S_{k-1}}-1}{2^{k-1}} \ge \frac{1}{2^k}.$$
\end{proof}

\newpage 

\section{Explicit Bounds and the Result in Dimension $7$} \label{Dim7}
	\indt Continue to suppose that $(R, \mf{m}_R, \kappa)$ is a complete normal characteristic $p > 2$ domain with algebraically closed residue field.  Assume that $R$ is not regular and that $\dim R = 7.$  The W-Y Conjecture is proven for complete intersections in \cite{enescu2005upper}, and therefore we assume that $R$ is not CI.  \\

For $R$ of dimension $7$ satisfying these conditions, the W-Y Conjecture reads:
\begin{align} \label{conjinc}
	\ehk{R} > \ehk{R_{p, 7}} \ge \frac{332}{315}.
\end{align}
A (tedious) computation following the techniques of \cite{HanMonskySomeSH} shows that 
\begin{align*}
	\ehk{R_{p, 7}} = \frac{332p^4+304p^2+192}{315p^4+273p^2+168} = \frac{332}{315} + \frac{244p^2+224}{4725p^4+4025p^2+2520},
\end{align*}
which is enough to establish the second inequality.  \\

Noting that $\ehk{R_{p, 7}}$ is a decreasing function\footnote{This can be seen immediately from the equality $$\frac{d}{dp} \ehk{R_{p, 7}} = - \frac{488 p^5+896 p^3+128 p}{4725 p^8+8190 p^6+8589 p^4+4368 p^2+1344}$$} of $p \ge 3$, we now use the techniques developed in the first part of this paper to establish the inequality,
\begin{align*}
	\ehk{R_{3, 7}} = 71/67 < \ehk{R},
\end{align*}
for all $R$ satisfying the conditions above.  This shows that (\ref{conjinc}) holds for all characteristics $p \ge 3$ and, as was shown in section \ref{reductionblurb}, is enough to settle the conjecture in dimension $7$ in general. \\

First we settle some preliminary cases.  Numerical computations and graphics were done using \texttrademark{MATHEMATICA}.
\subsection{Special Values of $e$ and $\mu$} 

\indt Since $R$ is singular, we have $e(R) \ge 2.$  Moreover, if the Hilbert-Samuel multiplicity of $R$ satisfies $$2 \le e(R) \le 5,$$ the conjecture holds by the results of \cite{IanNewEst}.\footnote{If $e(R)=2$, then, by Theorem 2.4(iii) of \cite{IanNewEst}, either $\ehk{R} \ge \frac{2}{2-1} = 2$ or $R$ is Cohen-Macauly, in which case it is a hypersurface and the result is known by \cite{enescu2005upper}.  For $e(R) = 3, 4$ and $5$, see the detailed discussion in \cite{IanNewEst} after the proof of Theorem 4.1.} \\

The inequality, $$\ehk{R} \ge \frac{e(R)}{d!},$$ is well known to hold for rings of dimension $d$, and therefore the conjecture is also known in $\dim R = 7$ if $$e(R) > \left(\frac{71}{67}\right)7! = 5340.9.$$  Hence, we need only worry about $R$ with Hilbert-Samuel multiplicity that satisfies,  $$6 \le  e(R) \le 5340.$$ \\ 

Now, for what follows, fix a minimal reduction, $\unideal{x} \subset \mf{m}_R,$ of the maximal ideal in $R$, and write $$\mu(R) = \mu_{R} \left(\mf{m}_R/\unideal{x}^*\right)$$ for the minimal number of generators of $\mf{m}_R/\unideal{x}^*.$  The argument at the beginning of the proof of Theorem 5.2 in \cite{IanNewEst} shows that we can assume\footnote{They show that either $R$ is Cohen-Macaulay with minimal multiplicity, in which case the conjecture is known (see Theorem 2.4(iv) of \cite{IanNewEst}), or the inequality $\mu(R) \le e(R) - 2,$ holds.}  $$\mu(R) \le e(R)-2.$$

If $\mu(R) = 1,  2$ or $3$, then the bound in Theorem 3.2 of \cite{IanNewEst}, gives, for $e(R) \ge 6,$\footnote{The bounds listed for $\mu(R) =1, 2$ and $3$ are achieved at $s = 4$, $s = 3.56745$, and $s = 3.32317$ respectively.} $$\ehk{R} \ge 6\left(\nu_{s} - \nu_{s-1}\right) \ge 2.87619 \text{   when } \mu(R) = 1,$$ $$\ehk{R} \ge 6\left(\nu_{s} - 2\nu_{s-1}\right) \ge 1.84215 \text{   when } \mu(R) = 2,$$ $$\text{and  } \ehk{R} \ge 6\left(\nu_{s} - 3\nu_{s-1}\right) \ge 1.33532 \text{   when } \mu(R) = 3,$$ and thus the conjecture is satisfied in these cases.
\subsection{New Bounds in Dimension $7$} \label{NewDim7Bnds}
\indt Now we turn our attention to the remaining cases.  Continuing with $(R, \mf{m}_R, \kappa)$ a complete normal characteristic $p > 2$ domain with algebraically closed residue field, we now consider the cases $$6 \le e(R) \le 5340,$$ and $$4 \le \mu(R) \le e(R)-2.$$ \\

Choose $z = z_1, z_2 \in \mf{m}_R \setminus \mf{m}_R^2$ which are part of a common minimal reduction of $\mf{m}_R$ and map to minimal generators of $\mf{m}_R/\unideal{x}^*.$\footnote{This is always possible as long as $\mu(R) \ge 2$ --- see Theorem 8.6.6 of \cite{HSintegral}.}  Letting $S_1 = R[z_2^{1/2}]$, recall that, by lemma \ref{NotNormBnd}, if $S_1$ is not normal, then $$\ehk{R} \ge 1 + \frac{1}{2},$$ so the conjecture is more than satisfied in this case. \\

Thus we can assume that $S_1$ is normal and take $k=1$ in Theorem \ref{BigBoundThm}.  Defining
\begin{align*}
	B\left(e, \mu, s, t\right) = 1 - \frac{t}{2} + e\left(\nu_s - (\mu-2)\nu_{s-1} - \nu_{s-\frac{1}{2}} - \nu_{s-t} \right),
\end{align*} 
our bound takes the form $$\ehk{R} \ge B\left(e(R), \mu(R), s, t\right)$$ for all $s \ge 0$ and $t\in [0,1].$ \\

Recall that we are concerned with $4 \le \mu(R) \le e(R)-2$ and note that, for fixed values of $e, s$ and $t$, the bounding function, $B$, increases when $\mu$ gets smaller.  With this observation in mind, define
\begin{align*}
	H_e\left(s, t\right) = B\left(e, e-2, s, t\right),
\end{align*}
so that we have
\begin{align*}
	\ehk{R} \ge B\left(e(R), \mu(R), s, t\right) \ge H_{e(R)}\left(s, t\right).
\end{align*}
  Our goal is now to show that, for each $6 \le e \le 5340$,  $$ \max \bigg ( H_e\left(s, t\right) \, \bigg | \, 0 \le s, \, \, 0 \le t \le 1 \bigg ) > 71/67 \approx 1.05970149.$$

The bounding function $H_e\left(s, t\right)$ is piece-wise polynomial in $s$ and $t$, but it is made up of a great many pieces.  Working with it by hand has proven to be extremely tedious.  Instead, we utilized the numerical algorithms of \texttrademark{MATHEMATICA} to find the optimal values of $H_e$ and the corresponding values of $s$ and $t$. 

\newpage

 \noindent\makebox[\linewidth]{\rule{6in}{0.4pt}}
				
		\begin{multicols}{2} 
		\hspace{.1in} \\ 
		
			   \vspace{.1in} 
				This image illustrates our bound in the $e = 7$ case.  A truncated graph of $H_7(s,t)$ appears in yellow, while the W-Y bound of $321/315$ is shown in blue. \\ 
				\vspace{.15in} 
				
				$H_7$ achieves the value $1.06056$
				at $s = 2.74118, \, t=0.779643$ \\
				\begin{flushright}
						\includegraphics[scale=.55]{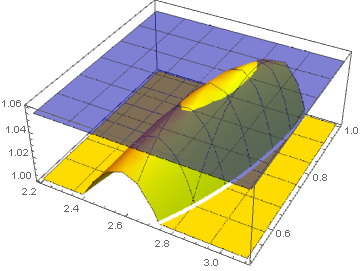}
				\end{flushright}
		\end{multicols}
		\begin{center}{Figure 1}\end{center}
	\begin{center}\rule{6in}{0.4pt}\end{center}

\indt	The table below lists the bounds achieved for $H_e(s,t)$ for $e = 6$ through $e = 12$.  \\

\begin{center}
	\begin{tabular}{|c|c|c|} 
		\hline 
		$e(R) $ & $(s, t) $ & $\ehk{R} \ge H_e(s, t)$ \\
		\hline
		\hline 
		$6$ & $(2.84243, 0.8)$ & $1.06447$ \\
		$7$ & $(2.74118, 0.779643)$ & $1.06056$ \\
		$8$ & $(2.65255, 0.739206)$ & $1.06024$ \\
		$9$ & $(2.58286, 0.710503)$ & $1.06183$ \\
		$10$ & $(2.52575, 0.688955)$ & $1.06438$ \\
		$11$ & $(2.47759, 0.672106)$ & $1.06742$ \\
		$12$ & $(2.43609, 0.658519)$ & $1.07073$ \\
		\hline  
		\hline 
	\end{tabular}
\end{center}

\vspace{.2in}

	\begin{center}{Table 1}\end{center}

\begin{center}\rule{6in}{0.4pt}\end{center}

\vspace{.25in}

To handle the remaining cases, we follow \cite{IanNewEst} in observing that the bounding function, 
\begin{align*}
	H_e(s,t) &= 1 - \frac{t}{2} + e\left(\nu_s - (e-4)\nu_{s-1} - \nu_{s-\frac{1}{2}} - \nu_{s-t} \right) \\
	 &= -e^2 \nu_{s-1} + e \left(\nu_{s} + 4\nu_{s-1} - \nu_{s-\frac{1}{2}} - \nu_{s-t}\right) + 1 - \frac{t}{2}
\end{align*}
is a downward facing parabola when viewed as a function of $e$. \\

In particular, for $(s, t) = (s_0, t_0)$ fixed, the bounding function obtains a max of $$H_{e_{\text{max}}}(s_0, t_0),$$ where $$e_{\text{max}} = \dfrac{\nu_{s_0}+4\nu_{s_0-1}-\nu_{s_0-\frac{1}{2}}-\nu_{s_0-t_0}}{2\nu_{s_0-1}}.$$   

Furthermore, given any $e_1$ and $e_2$ with $$e_1 \le e_{\text{max}} \le e_2$$ if $R$ is any ring under consideration such that $$e_1 \le e(R) \le e_2,$$ then we are able to conclude that $$\ehk{R} \ge H_{e(R)}(s_0, t_0) \ge \min \left(H_{e_1}(s_0, t_0), \, H_{e_2}(s_0, t_0)\right).$$  We use this to bound a range of multiple $e$-values at a time.  The results are recorded in the following table.\footnote{The bounds in this table are not the best that our methods achieve for each specific value of $e \ge 13$.  The particular values of $e_1$ and $e_2$, and the corresponding bounds on $e_{HK}$, were chosen, in an 'informed trial and error' process, to make the overall presentation of the data as palatable as possible.}

\vspace{.25in}

\begin{center}\rule{6in}{0.4pt}\end{center}

\indt	The table below lists the bounds achieved for $H_e(s,t)$ for ranges of $e$-values between $e = 13$ and $e = 5340$.  \\

This completes the proof of Theorem \ref{WYDim7}. $\qed$

\begin{center}
	\begin{tabular}{|l|c|c|c|} 
		\hline 
		$[e_1, e_2] $ & $(s_0, t_0) $ & $e_{\text{max}}$ & $e_{HK} \ge \min \left( H_{e_1}, \, H_{e_2} \right)$ \\
		\hline
		\hline 
		$[13, 19]$ & $(2.34, 0.62)$ & $15.973$ & $1.06843$ \\
		$[20, 40]$ & $(2.12, 0.6)$ & $31.2399$ & $1.07266$ \\
		$[41, 105]$ & $(1.9, 0.55)$ & $72.3972$ & $1.12153$ \\
		$[106, 227]$ & $(1.75, 0.5)$ & $151.062$ & $1.20165$ \\
		$[228, 650]$ & $(1.6, .475)$ & $402.416$ & $1.32149$ \\
		$[651, 1600]$ & $(1.5, 0.45)$ & $937.946$ & $1.45925$ \\
		$[1601, 5340]$ & $(1.375, 0.41)$ & $3891.82$ & $2.84311$ \\
		\hline  
		\hline 
	\end{tabular}
\end{center}

\vspace{.2in}

	\begin{center}{Table 2}\end{center}

\begin{center}\rule{6in}{0.4pt}\end{center}
\vspace{.25in}

\newpage 

\section{Conclusion}
\subsection{Higher Dimensions}
\indt The arguments in this paper work in any dimension $\ge 2.$  However, the bound in Theorem \ref{BigBoundThm}, which is our primary weapon, quickly looses potency as the dimension of $R$ increases.  We have only been able to get partial results in dimensions $8$ and above. \\

For example, taking $k=4$ in Theorem \ref{BigBoundThm} and arguing as we did for the dimension $7$ case, our bounds appear to be strong enough to verify the conjecture in dimension $8$ as long as $e(R) \ge 21.$  Unfortunately, when $R$ has dimension $8,$ and $6 \le e(R) \le 20,$ our tools have so far only been strong enough to establish weaker forms of the conjecture --- our best cursory calculations (with various values for $k$) attain the desired bounds in dimension $8$ when either $e(R) \ge 21,$ as mentioned, or $\mu(R) \le e-6,$ though we will not go into detail about this here. \\

Similar computations indicate that our methods can establish the conjecture in dimension $9$ for $e(R) \ge 49,$ in dimension $10$ for $e(R) \ge 135,$ and in dimension $11$ for $e(R) \ge 632.$

\noindent\makebox[\linewidth]{\rule{6in}{0.4pt}}
				
		\begin{multicols}{2} 
		\hspace{.1in} \\ 
				
			   \vspace{.1in} 
				This image illustrates our bound for $e = 21,$ $\mu = 19,$ in dimension $d=8,$ with $k=4$.  A truncated graph of the bounding function appears in yellow, while the characteristic-free dimension $8$ W-Y bound of $8341/8064 \approx 1.03435$ is shown in blue. \\ 
				\vspace{.15in} 
				
				The bounding function has value $1.03545$
				at $s = 2.17991, \, t=0.706957$ \\
				\begin{flushright}
						\includegraphics[scale=.55]{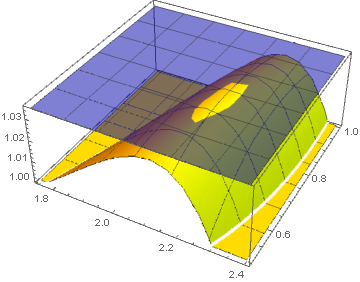}
				\end{flushright}
		\end{multicols}
		\begin{center}{Figure 2}\end{center}	
	\begin{center}\rule{6in}{0.4pt}\end{center}
	\vspace{.1in}

\newpage

\subsection{The $\varphi_J$ Function}

\indt The function, $\PsifunZ{t}$, studied in section \ref{PhiSection} is very interesting and there are good reasons to expect that a more detailed understanding of this function may lead to improvements to the bounds studied in this paper.  We have been able to explicitly compute values of $\PsifunZ{t}$ for some simple diagonal hypersurfaces using the 'Representation Ring' methods introduced in \cite{HanMonskySomeSH, MonskyTeix2006p} and hope to gain a better understanding of $\varphi$ and it's derivatives from these examples.

\subsection{Acknowledgments} 
\indt The authors are grateful to the University of Tulsa and the University of Missouri, Columbia.  They would also like to thank \texttrademark{ZOOM} for making this collaboration possible without extensive commuting.

\newpage

\bibliography{bib}

\newcommand{\etalchar}[1]{$^{#1}$}
\begin{thebibliography}{CDHZ12}

\bibitem[AE08]{aberbach2008lower}
Ian~M Aberbach and Florian Enescu.
\newblock Lower bounds for {Hilbert-Kunz} multiplicities in local rings of
  fixed dimension.
\newblock {\em Michigan Mathematical Journal}, 57:1--16, 2008.

\bibitem[AE13]{IanNewEst}
Ian~M. Aberbach and Florian Enescu.
\newblock {New estimates of {Hilbert–Kunz} multiplicities for local rings of
  fixed dimension}.
\newblock {\em Nagoya Mathematical Journal}, 212(none):59 -- 85, 2013.

\bibitem[BST12]{blickle2012f}
Manuel Blickle, Karl Schwede, and Kevin Tucker.
\newblock F-signature of pairs and the asymptotic behavior of {F}robenius
  splittings.
\newblock {\em Advances in Mathematics}, 231(6):3232--3258, 2012.

\bibitem[CDHZ12]{CHDZ}
Olgur Celikbas, Hailong Dao, Craig Huneke, and Yi~Zhang.
\newblock Bounds on the {Hilbert-Kunz} multiplicity.
\newblock {\em Nagoya Mathematical Journal}, 205:149--165, 2012.

\bibitem[CK20]{cho2020volume}
Yunhi Cho and Seonhwa Kim.
\newblock Volume of hypercubes clipped by hyperplanes and combinatorial
  identities.
\newblock {\em The Electronic Journal of Linear Algebra}, 36:228--255, 2020.

\bibitem[CS21]{CS20}
Nick Cox-Steib.
\newblock $\mathfrak{m}$-adic perturbations in {N}oetherian local rings.
\newblock {\em Thesis, University of Missouri}, 2021.

\bibitem[ES05]{enescu2005upper}
Florian Enescu and Kazuma Shimomoto.
\newblock On the upper semi-continuity of the {Hilbert--Kunz} multiplicity.
\newblock {\em Journal of Algebra}, 285(1):222--237, 2005.

\bibitem[GM10]{gessel2010limit}
Ira~M Gessel and Paul Monsky.
\newblock {The limit as $p\to \infty$ of the {H}ilbert-{K}unz multiplicity of
  $\Sigma x_i^{d_i}$}.
\newblock {\em arXiv preprint:1007.2004}, 2010.

\bibitem[HM93]{HanMonskySomeSH}
C.~Han and Paul Monsky.
\newblock Some surprising {H}ilbert-{K}unz functions.
\newblock {\em Mathematische Zeitschrift}, 214:119--135, 1993.

\bibitem[HSPS06]{HSintegral}
C.~Huneke, I.~Swanson, Cambridge~University Press, and London~Mathematical
  Society.
\newblock {\em Integral Closure of Ideals, Rings, and Modules}.
\newblock Cambridge University Press, New York, 2006.

\bibitem[Hun14]{huneke2014hilbert}
Craig Huneke.
\newblock Hilbert-{K}unz multiplicities and the {F}-signature.
\newblock {\em arXiv preprint arXiv:1409.0467}, 2014.

\bibitem[HW{\etalchar{+}}20]{hug2020lectures}
Daniel Hug, Wolfgang Weil, et~al.
\newblock {\em Lectures on convex geometry}.
\newblock Springer, Switzerland, 2020.

\bibitem[iWiY05]{WY3Dim}
Kei ichi Watanabe and Ken ichi Yoshida.
\newblock Hilbert-{K}unz multiplicity of three-dimensional local rings.
\newblock {\em Nagoya Mathematical Journal}, 177:47 -- 75, 2005.

\bibitem[JNS{\etalchar{+}}23]{jeffries2023lower}
Jack Jeffries, Yusuke Nakajima, Ilya Smirnov, Kei-Ichi Watanabe, and Ken-Ichi
  Yoshida.
\newblock Lower bounds on {Hilbert--Kunz} multiplicities and maximal
  {F}-signatures.
\newblock In {\em Mathematical Proceedings of the Cambridge Philosophical
  Society}, volume 174, pages 247--271. Cambridge University Press, 2023.

\bibitem[Kun69]{kunz1969characterizations}
Ernst Kunz.
\newblock Characterizations of regular local rings of characteristic p.
\newblock {\em American Journal of Mathematics}, 91(3):772--784, 1969.

\bibitem[Kun76]{kunz1976noetherian}
Ernst Kunz.
\newblock On {N}oetherian rings of characteristic p.
\newblock {\em American Journal of Mathematics}, 98(4):999--1013, 1976.

\bibitem[MM08]{marichal2006slices}
Jean-Luc Marichal and Michael~J Mossinghoff.
\newblock Slices, slabs, and sections of the unit hypercube.
\newblock {\em Online Journal of Analytic Combinatorics}, 3:1--11, 2008.

\bibitem[Mon83]{monsky1983hilbert}
Paul Monsky.
\newblock The {H}ilbert-{K}unz function.
\newblock {\em Mathematische Annalen}, 263(1):43--49, 1983.

\bibitem[MT06]{MonskyTeix2006p}
Paul Monsky and Pedro Teixeira.
\newblock p-fractals and power series—ii.: Some applications to
  {Hilbert--Kunz} theory.
\newblock {\em Journal of algebra}, 304(1):237--255, 2006.

\bibitem[PS18]{polstra_smirnov_2018}
Thomas Polstra and Ilya Smirnov.
\newblock Continuity of {H}ilbert–{K}unz multiplicity and {F}-signature.
\newblock {\em Nagoya Mathematical Journal}, 239:1–--24, 2018.

\bibitem[Ser12]{serrelocal}
Jean-Pierre Serre.
\newblock {\em Local algebra}.
\newblock Springer, New York, 2012.

\bibitem[Tri21]{trivedi2021lower}
Vijaylaxmi Trivedi.
\newblock The lower bound on the {HK} multiplicities of quadric hypersurfaces.
\newblock {\em arXiv preprint arXiv:2101.03905}, 2021.

\bibitem[Tri23]{trivedi2023hilbert}
Vijaylaxmi Trivedi.
\newblock The {Hilbert-Kunz} density functions of quadric hypersurfaces.
\newblock {\em Advances in Mathematics}, 430:109--207, 2023.

\bibitem[WY00]{watanabe2000hilbert}
Kei-ichi Watanabe and Ken-ichi Yoshida.
\newblock Hilbert--{K}unz multiplicity and an inequality between multiplicity
  and colength.
\newblock {\em Journal of Algebra}, 230(1):295--317, 2000.

\bibitem[Yos09]{yoshida2009small}
K-i Yoshida.
\newblock Small {Hilbert-Kunz} multiplicity and {(A1)}-type singularity.
\newblock In {\em Proceedings of the 4th Japan-Vietnam Joint Seminar on
  Commutative Algebra by and for Young Mathematicians}, 2009.

\end{thebibliography}
\bibliographystyle{alpha}

\end{document}